# RESOLVABLE DESIGNS WITH LARGE BLOCKS

By J. P. Morgan[1] and Brian H. Reck[2]

*Virginia Tech and University of Pittsburgh*

*Dedicated to the memory of I. M. Chakravarti*

Resolvable designs with two blocks per replicate are studied from an optimality perspective. Because in practice the number of replicates is typically less than the number of treatments, arguments can be based on the dual of the information matrix and consequently given in terms of block concurrences. Equalizing block concurrences for given block sizes is often, but not always, the best strategy. Sufficient conditions are established for various strong optimalities and a detailed study of E-optimality is offered, including a characterization of the E-optimal class. Optimal designs are found to correspond to balanced arrays and an affine-like generalization.

**1. Introduction.** Block designs arise in comparative experimentation as fundamental devices for improving efficiency when working with heterogeneous experimental units. The blocks are simply a partition of the units into (say) $b$ sets exhibiting homogeneity within sets. Given the blocks, of sizes $k_1, k_2, \ldots, k_b$, a *block design* is an assignment of $v$ treatments to the $\sum_{j=1}^{b} k_j$ units. Optimality theory for block designs attempts to determine which of the many possible assignments is in some sense best.

In some applications there are restrictions on the collection of possible assignments. A block design is *resolvable* if the blocks can be partitioned into *replicates*, defined as sets of blocks with the property that each treatment is assigned to one unit in each set. The practical impact of, and motivation for, resolvability is to gain orthogonality between treatments and nuisance factors of concern. For instance, resolvability in sequential experimentation,

---

Received February 2005; revised February 2006.
[1]Supported by NSF Grant DMS-01-04195.
[2]Supported by NIMH Training Grant 5T32MH020053-05.

with replicates corresponding to time periods, is used to mitigate time effects. Resolvability can likewise be useful in multi-site experiments and in experiments with multiple individuals handling experimental runs. Notably, the United Kingdom has for some time required the use of resolvable designs in agricultural field trials (see [20]).

The combinatorial study of resolvability in block designs goes back at least as far as the well-known Kirkman's [14] schoolgirl problem. The notion entered the statistical lexicon with Yates' work on square *lattice* designs [29, 30], although the term "resolvable design" was introduced by Bose [4]. Yates' lattice designs were extended to rectangular lattices by Harshbarger [11, 12]; see also [2]. Williams [27] and Patterson and Williams [21] introduced a large family of resolvable designs which they termed $\alpha$-designs. Williams, Patterson and John [28] derived resolvable designs with two replicates from BIBDs (balanced incomplete block designs). Bailey, Monod and Morgan [1] proved strong optimality for the affine resolvable designs introduced by Bose [4]. Resolvable BIBDs have received significant attention in both the combinatorial and the statistical literature; for a summary, see Morgan [17], who surveys the major classes of resolvable designs with many references.

With $\alpha$-designs, Patterson and Williams [21] provided a flexible method for obtaining reasonably efficient resolvable designs for a wide range of values $v$, constant block size $k$ and replication $r$. They also adapted their method to obtain resolvable designs with two different block sizes, $k$ and $k-1$, a first attempt at addressing the obvious restriction that resolvability with equal block sizes can be achieved only for $v$ a multiple of $k$. John, Russell, Williams and Whitaker [13], in revisiting that idea, concluded that the $\alpha$-technique for two block sizes could produce relatively inefficient designs for small $v$ and recommended an interchange algorithm for construction of designs with better efficiency. John et al. [13] also discussed the practical need for resolvable designs with unequal block sizes; for example, about half of 245 experiments examined by Patterson and Hunter [19] had unequal block sizes. See also [20].

To the authors' knowledge, the literature contains no systematic work on determining optimal resolvable block designs when block sizes need not be equal. This paper will undertake such work, for the special case of two blocks per replicate. Only with two blocks in each replicate must a block be large in the sense of containing at least half of the treatments. And with two blocks per replicate, the block sizes must be unequal for any odd $v$ (although this work is not restricted to odd $v$).

Let $\mathcal{D}(v, r; k_1, k_2)$ denote the class of all resolvable block designs with $r$ replicates of $v$ treatments, each replicate consisting of two blocks, of sizes $k_1$ and $k_2$, where $k_1 + k_2 = v$. Of special interest is $k_1 = k$, $k_2 = k - 1$ for odd $v = 2k - 1$, but here no restrictions are placed on the two block sizes for the general theory. Let $k_1$ denote the larger block size, $k_1 \geq k_2$, so that with no loss of generality $\frac{v}{2} \leq k_1 \leq v-2$. In the most common applications of resolvable designs, the number of treatments is large relative to the number of replicates; here, $r \leq v-1$ is required, allowing optimality problems to be more easily attacked using the dual of the information matrix. These considerations frame the goal of this paper: to determine the best design $d \in \mathcal{D}(v,r;k_1,k_2)$. An example of a resolvable design in $\mathcal{D}(9,5;5,4)$ is shown in Table 1, with the blocks written as columns. Later, this design will be proven optimal with respect to many useful criteria.

TABLE 1
*A resolvable design in $\mathcal{D}(9,5;5,4)$*

| | | | | |
|---|---|---|---|---|
| 1 6 | 1 2 | 1 4 | 1 3 | 1 2 |
| 2 7 | 4 3 | 2 5 | 2 4 | 3 5 |
| 3 8 | 5 6 | 3 6 | 5 6 | 4 6 |
| 4 9 | 7 8 | 7 8 | 8 7 | 7 9 |
| 5 | 9 | 9 | 9 | 8 |

For any resolvable design $d \in \mathcal{D}$, ignoring the replicate grouping leaves an underlying *simple block design* for $v$ treatments in $2r$ blocks. If the roles of blocks and treatments are reversed in this underlying design, another simple block design with $2r$ treatments in $v$ blocks of size $r$ is produced, called the *dual design*.

Section 2 provides optimality background and establishes equivalences with the dual problem. Section 3 finds conditions for global optimality of designs having equality of block concurrences; for those who wish to jump ahead, the main results there are Theorems 4–7. Section 4 characterizes the E-optimal designs (Theorem 14) and finds the Schur-best of the E-optimal designs (Theorem 15). The main results are applied to special cases with $k_1 - k_2 \leq 2$ in Section 5 and designs for these cases are constructed in Section 6. Concluding comments appear in Section 7.

**2. Model, information and optimality criteria.** Let $y_{hjl}$ denote the yield from the $l$th experimental unit in block $j$ of replicate $h$. Thus, the triples $(h,j,l)$ identify the experimental units and the design $d$ corresponds to a map $d[h,j,l]$ from the units to the set of treatments. The standard linear model for the yields incorporates a mean effect $\mu$, replicate effects $\rho_h$, block effects $\beta_{hj}$, treatment effects $\tau_{d[h,j,l]}$ and mean zero, uncorrelated, equivariable random error terms $e_{hjl}$:

$$y_{hjl} = \mu + \rho_h + \beta_{hj} + \tau_{d[h,j,l]} + e_{hjl},$$

$h = 1, \ldots, r$; $j = 1, 2$; $l = 1, \ldots, k_j$. This model may be written in matrix terms as

$$y = \mu 1 + (I_r \otimes 1_v)\rho + L\beta + X_d\tau + e,$$

where, with the $y_{hjl}$ lexicographically ordered in the $vr \times 1$ yield vector $y$, the block incidence matrix is $L_{vr \times 2r} = I_r \otimes \begin{pmatrix} 1_{k_1} & 0_{k_1} \\ 0_{k_2} & 1_{k_2} \end{pmatrix}$ and the design matrix is the $vr \times v$ incidence matrix $X_d$, for which row $(h, j, l)$ has a 1 in column $i$ if and only if $d[h, j, l] = i$ [i.e., unit $(h, j, l)$ receives treatment $i$] and all other entries are zero. Replicate effects, block effects, treatment effects and error vectors are $\rho_{r \times 1}$, $\beta_{2r \times 1}$, $\tau_{v \times 1}$ and $e_{vr \times 1}$, respectively. Choice of design is equivalently choice of $X_d$.

Linear model theory says that the information matrix $C_d$ for estimation of the treatment effects $\tau$ is

$$(1) \qquad C_d = X_d'[I - L(L'L)^{-1}L']X_d = rI - N_d D_s^{-1} N_d',$$

where $D_s = L'L = \text{Diag}(k_1, k_2, k_1, k_2, \ldots, k_1, k_2)$ is the diagonal matrix of block sizes. The $v \times 2r$ matrix $N_d$ is the treament/block incidence matrix. The general off-diagonal element $(N_d N_d')_{i,i'}$ of $N_d N_d'$ is the number of blocks to which both treatments $i$ and $i'$ are assigned, called a *treatment concurrence count*. Notice that the replicate incidence matrix $(I_r \otimes 1_v)$ plays no role in (1); the same form of information matrix is obtained for any simple block design. It follows immediately that the information matrix for the corresponding dual design is

$$(2) \qquad C_{\text{dual}} = D_s - \frac{1}{r} N_d' N_d.$$

The off-diagonal elements of $N_d' N_d$ are *block concurrence counts*.

All treatment contrasts are estimable with design $d$ if and only if $C_d$ has rank $v - 1$. Any such $d$ is said to be *connected*; only connected designs are considered here. Most (but not all) commonly employed optimality criteria, including those to be used here, are functions of the $v - 1$ nonzero eigenvalues of $C_d$. These will be ordered and labeled $z_{d1} \leq z_{d2} \leq \cdots \leq z_{d,v-1}$.

Much of the optimality work below focuses on minimizing functions of the form

$$(3) \qquad \psi_f(z_d) = \sum_{i=1}^{v-1} f(z_{di}),$$

where $f$ is convex and $z_d$ is the vector of nonzero eigenvalues. If $f$ in (3) is Schur-convex (Schur-convex functions include the convex functions; see [3], Section II.3), then (3) is said to be a *Schur-criterion*. Design $d_1$ is *Schur-better* than design $d_2$ ($d_2$ is *Schur-inferior* to $d_1$) if $\psi_f(z_{d_1}) \leq \psi_f(z_{d_2})$ for all Schur-convex $f$ with strict inequality for at least one such $f$. A design optimal with respect to all (i.e., minimizing all) Schur-convex criteria is said to be *Schur-optimal*. If $f$ in (3) satisfies (i) $f$ is continuously differentiable on $(0, \max_{d \in \mathcal{D}} \text{tr}(C_d))$ with $f' < 0$, $f'' > 0$ and $f''' < 0$ and (ii) $\lim_{x \to 0} f(x) = \infty$, then (3) is said to be a *type-1 criterion* (see [7]). A design optimal with respect to all type-1 criteria is said to be *type-1-optimal*. One popular criterion belonging to both families just defined is the *A-criterion* specified by $f(x) = \frac{1}{x}$. A criterion not of the form (3) (though it can be written as a limit of such criteria) is

$$\psi_E(z_d) = \frac{1}{z_{d1}}, \tag{4}$$

called the *E-criterion*. An E-optimal design minimizes (4) or, equivalently, maximizes $z_{d1}$. For a broader discussion of optimality criteria and their statistical meanings, see [24].

For the current endeavor, it is advantageous to approach the $z_{di}$, and consequently (3) and (4), through $C_{\text{dual}}$. How this is done is shown in Lemma 1 below. Let $\phi_{dhh'}$ be the block concurrence for (i.e., the number of treatments common to) the blocks of size $k_1$ in replicates $h$ and $h'$. For $h \neq h'$, let $\phi^*_{dhh'} = \phi_{dhh'} - \frac{k_1^2}{k_1+k_2}$ and define the (symmetric) *optimality matrix* $M_d$ by

$$M_d = \begin{pmatrix} \frac{k_1 k_2}{v} & \phi^*_{d12} & \phi^*_{d13} & \cdots & \phi^*_{d1r} \\ & \frac{k_1 k_2}{v} & \phi^*_{d23} & \cdots & \phi^*_{d2r} \\ & & \frac{k_1 k_2}{v} & \cdots & \phi^*_{d3r} \\ & & & \ddots & \vdots \\ & & & & \frac{k_1 k_2}{v} \end{pmatrix}, \tag{5}$$

with eigenvalues $e_{d1} \geq e_{d2} \geq \cdots \geq e_{dr}$.

LEMMA 1. *For any $d \in \mathcal{D}(v, r; k_1, k_2)$, the eigenvalues of $C_d$ are $0$, $\max\{0, v - r - 1\}$ copies of $r$ and $(r - \frac{ve_{d1}}{k_1 k_2}, r - \frac{ve_{d2}}{k_1 k_2}, \ldots, r - \frac{ve_{dr}}{k_1 k_2})$. Consequently, a $\psi_f$-optimal design minimizing (3) will equivalently minimize $\sum_{h=1}^r f(r - \frac{ve_{dh}}{k_1 k_2})$ and an E-optimal design minimizing (4) will equivalently minimize $e_{d1}$.*

Lemma 1 is proved in Appendix A.1. Provided $r \leq v - 1$ (as earlier required in Section 1), $C_d$ has $v - r - 1$ eigenvalues fixed at $r$. Working with the dual through $M_d$ not only makes this evident, but allows these structurally fixed eigenvalues to be easily set aside.

With the optimality problem recast in terms of $M_d$ and its eigenvalues, a crucial concept for proofs of Schur-optimality is now defined.

DEFINITION. Let $\{x_i\}_{i=1}^n$ and $\{y_i\}_{i=1}^n$ be nonincreasing sequences of real numbers such that $\sum_{i=1}^n x_i = \sum_{i=1}^n y_i$. If

$$\sum_{i=1}^l x_i \leq \sum_{i=1}^l y_i \quad \text{for all } 1 \leq l \leq n$$

or, equivalently,

$$\sum_{i=n}^{n-l+1} x_i \geq \sum_{i=n}^{n-l+1} y_i \quad \text{for all } 1 \leq l \leq n,$$

then $\{y_i\}_{i=1}^n$ is said to *majorize* $\{x_i\}_{i=1}^n$.

The importance of majorization is evident in the following result (see, e.g., [3], page 40).

THEOREM 2. *Let $\{x_i\}_{i=1}^n$ and $\{y_i\}_{i=1}^n$ be nonincreasing sequences of real numbers such that $\sum_{i=1}^n x_i = \sum_{i=1}^n y_i$. Then $\sum_{i=1}^n f(x_i) \leq \sum_{i=1}^n f(y_i)$ for all real-valued convex functions $f$ if and only if $\{y_i\}_{i=1}^n$ majorizes $\{x_i\}_{i=1}^n$.*

If the sequences $\{x_i\}_{i=1}^n$ and $\{y_i\}_{i=1}^n$ are written as the elements of vectors $x$ and $y$, then the fact that $y$ majorizes $x$ is written as $y \succ x$ or $x \prec y$.

Let $e_d$ and $e_{d^*}$ be the vectors of eigenvalues for the optimality matrices for designs $d$ and $d^*$, respectively. Lemma 3 lists simple majorization facts ([3], page 30) used in subsequent sections. For majorization comparisons, the third of these allows work directly with the $e_{dh}$ rather than the $r - \frac{ve_{dh}}{k_1 k_2}$ (see Lemma 1). Thus, if $e_d \succ e_{d^*}$ and the two vectors are not identical, then $d^*$ is Schur-better than $d$. Design $d^*$ is Schur-optimal if $e_d \succ e_{d^*}$ for every $d \in \mathcal{D}$.

LEMMA 3. *For real numbers $\{x_i\}_{i=1}^n$ and $\{y_i\}_{i=1}^n$ with $\sum_{i=1}^n x_i = \sum_{i=1}^n y_i$,*

  (i) *if $x_1 \geq x_2 = x_3 = \cdots = x_n$ and $y_1 \geq x_1$, then $\{y_i\}_{i=1}^n$ majorizes $\{x_i\}_{i=1}^n$;*
  (ii) *if $x_1 = x_2 = \cdots = x_{n-1} \geq x_n$ and $x_n \geq y_n$, then $\{y_i\}_{i=1}^n$ majorizes $\{x_i\}_{i=1}^n$;*
  (iii) *if $\{y_i\}_{i=1}^n$ majorizes $\{x_i\}_{i=1}^n$, then $\{a - \frac{y_i}{b}\}_{i=1}^n$ majorizes $\{a - \frac{x_i}{b}\}_{i=1}^n$ for any real $a, b$.*

**3. Equal concurrence designs and global optimality.** Among simple block designs, the BIBDs are Schur-optimal, a result which follows from equality of treatment concurrences inducing complete symmetry of the information matrix. The analogous notion for duals is equality of block concurrences, which this section explores for utility with resolvable designs. A resolvable design $d \in \mathcal{D}(v, r; k_1, k_2)$ having block concurrence counts $\phi_{d12} = \phi_{d13} = \phi_{d23} = \cdots = \phi_{d,r-1,r} = \theta$ for some $k_1 - k_2 \leq \theta \leq k_1$ is called an *equal concurrence design with common concurrence* $\theta$, or ECD($\theta$). The subsections below will explore, in turn, (i) ECDs with equality of eigenvalues in the optimality matrix, (ii) other ECDs which can be proven Schur-optimal and (iii) Schur-domination of designs with unequal block concurrences by one or more ECDs.

3.1. *Schur-optimality via equality of eigenvalues.* For an ECD($\theta$), the optimality matrix (5) is

$$(6) \qquad M_d = \left[\frac{p}{v} - \left(\theta - \frac{k_1^2}{v}\right)\right] I + \left(\theta - \frac{k_1^2}{v}\right) J,$$

where $I$ is the $r \times r$ identity matrix, $J$ is the $r \times r$ matrix of ones and $p = k_1 k_2$ is the product of the block sizes. Like a BIBD information matrix, it is completely symmetric, but unlike that matrix, $M_d$ for an ECD($\theta$) is nonsingular and thus can have two distinct, relevant eigenvalues, rather than just one.

THEOREM 4. *Suppose that $\mathcal{D}(v, r; k_1, k_2)$ is a resolvable design setting for which $(k_1 + k_2) | k_1^2$ and define*

$$(7) \qquad \theta^* = \frac{k_1^2}{k_1 + k_2} = \frac{k_1^2}{v}.$$

*Then ECD($\theta^*$)'s in $\mathcal{D}$ are Schur-optimal whenever they exist.*

PROOF. Given the conditions on $\mathcal{D}$, the inequalities $k_1 - k_2 \leq \frac{k_1^2}{k_1+k_2} \leq k_1$ imply that $\theta = \theta^*$ is an admissible value for the common block concurrence of an ECD($\theta$). For $\theta = \theta^*$, $M_d$ in (6) is $\frac{p}{v} I$. Since the eigenvalues of $M_d$ are all identical, they are majorized by the eigenvalues of every competing design. □

Corollary 3.4 of [1] established that affine-resolvable designs are Schur-optimal. Theorem 4 generalizes that result when there are two blocks per replicate. Here, the optimality condition is that all pairs of blocks of size $k_1$ have the same concurrence (7). When $k_1 = k_2$ and $2|k_1$, ECD($\theta^*$)'s are affine-resolvable designs. It is obvious that (6) has two distinct eigenvalues for any $\theta \neq \theta^*$.

EXAMPLE. Consider the setting $\mathcal{D}(9, 4; 6, 3)$. Then $(k_1 + k_2)|k_1^2$ with $\theta^* = 4$ and if an ECD(4) exists, it is Schur-optimal. In fact, an ECD(4) does exist and is shown in Table 2.

The settings for which $k_1^2$ is a multiple of $k_1 + k_2$ are relatively sparse (a situation much like that of BIBDs relative to all simple block design settings). For the 1225 pairs $2 \leq k_2 < k_1 \leq 51$, only 23 meet the divisibility requirement implied by (7). Theorem 4 is thus only a start, albeit an important one.

TABLE 2
*A Schur-optimal ECD(4)
in $\mathcal{D}(9,4;6,3)$*

| | | | |
|---|---|---|---|
| 1 7 | 1 5 | 1 3 | 3 1 |
| 2 8 | 2 6 | 2 4 | 4 2 |
| 3 9 | 3 9 | 5 9 | 5 9 |
| 4   | 4   | 6   | 6   |
| 5   | 7   | 7   | 7   |
| 6   | 8   | 8   | 8   |

3.2. *Global optimality of other ECDs.* Good designs are expected to have eigenvalue structures "close" to that of ECD($\theta^*$)'s, suggesting this question: Is some equal concurrence design Schur-optimal when $(k_1 + k_2) \nmid k_1^2$? To investigate this question, define the *block concurrence parameter* $\bar{\theta}$ by

$$(8) \qquad \bar{\theta} = \text{int}\left(\frac{k_1^2}{k_1 + k_2}\right)$$

and write

$$(9) \qquad \gamma = \frac{k_1^2}{k_1 + k_2} - \bar{\theta}.$$

Then $0 \leq \gamma < 1$ and a necessary condition for existence of ECD($\theta^*$) is $\gamma = 0$. Consequently, $\gamma$ is called the *block discordancy coefficient*; it measures the departure of the block sizes from that required for equality of all eigenvalues. The parameter $\gamma$ will play a pivotal role in the remainder of this paper, as will the concurrence discrepancies defined next.

Define the block concurrence *discrepancy matrix* $\Delta_d = (\delta_{dhh'})$, where

$$\delta_{dhh'} = \begin{cases} \phi_{dhh'} - \bar{\theta}, & \text{if } h \neq h', \\ 0, & \text{if } h = h'. \end{cases}$$

For $h \neq h'$, the off-diagonal elements $\delta_{dhh'}$ will be referred to as block concurrence *discrepancies*. Rewritten in terms of block discrepancies and the discordancy coefficient, the general optimality matrix (5) becomes

$$(10) \qquad M_d = \frac{p}{v} I - \gamma(J - I) + \Delta_d.$$

For ECD($\theta^*$), $\gamma = 0$, $\Delta_d = 0$ and $M_d = \frac{p}{v} I$. If $\gamma \neq 0$, (10) shows that the form of optimal design may depend both on the magnitude of $\gamma$ and the values of the discrepancies. Matrix (10) is completely symmetric for, and only for, ECDs, in which case $\Delta_d$ is an integer multiple of $J - I$.

The ECDs which are combinatorially closest to $\text{ECD}(\theta^*)$ are those with either $\theta = \bar{\theta}$ or $\theta = \bar{\theta} + 1$. Not surprisingly, these are strong competitors in the optimality race. $\text{ECD}(\bar{\theta})$s have $\phi_{dhh'} = \bar{\theta}$ for each $1 \leq h \neq h' \leq r$, so $\Delta_d = 0$ and the eigenvalues of $M_d$ are

$$\xi_1(\gamma) = \frac{p}{v} + \gamma \quad \text{and} \quad \xi_2(\gamma) = \frac{p}{v} - (r-1)\gamma, \tag{11}$$

with frequencies $r - 1$ and 1, respectively. The two eigenvalues satisfy $\xi_1(\gamma) \geq \xi_2(\gamma)$. If $\phi_{dhh'} = \bar{\theta} + 1$ for every $h \neq h'$, then $\text{ECD}(\bar{\theta} + 1)$'s have $\Delta_d = (J - I)$ and the eigenvalues of $M_d$ are

$$\begin{aligned}\xi_1(\gamma - 1) &= \frac{p}{v} - (1 - \gamma) \quad \text{and} \\ \xi_2(\gamma - 1) &= \frac{p}{v} + (r-1)(1-\gamma),\end{aligned} \tag{12}$$

with frequencies $r - 1$ and 1, respectively, and with $\xi_2(\gamma - 1) \geq \xi_1(\gamma - 1)$.

THEOREM 5. *For $0 \leq \gamma \leq \frac{1}{2}$, $\text{ECD}(\bar{\theta})$'s are type-1-optimal and E-optimal.*

PROOF. The eigenvalues of the information matrix for any design $d \in \mathcal{D}(v, r; k_1, k_2)$ are $0 < z_{d1} \leq z_{d2} \leq \cdots \leq z_{dr}$ and $v - r - 1$ copies of $r$, and $\sum_{i=1}^{r} z_{di} = r(r-1)$ is constant for all designs in $\mathcal{D}$. For an $\text{ECD}(\bar{\theta})$, call it $\bar{d}$, the $z_{\bar{d}i}$, following from (11), have the form $z_{\bar{d}1} = z_{\bar{d}2} = \cdots = z_{\bar{d},v-2} \leq z_{\bar{d},v-1}$. Theorem 2.3 of [7] thus gives the result if it can also be shown that $\bar{d}$ minimizes $\sum_{i=1}^{v-1} z_{di}^2$ over $\mathcal{D}$.

For $d \in \mathcal{D}(v, r; k_1, k_2)$ with optimality matrix $(M_d)_{hh'} = (\delta_{dhh'} - \gamma)$ having trace $\operatorname{tr} M_d = \frac{pr}{v}$,

$$\begin{aligned}\operatorname{tr} C_d^2 &= \sum_{i=1}^{v-1} z_{di}^2 = (v - r - 1)r^2 + \sum_{h=1}^{r}\left(r - \frac{ve_{dh}}{p}\right)^2 \\ &= (v-1)r^2 - \frac{2vr}{p} \operatorname{tr} M_d + \frac{v^2}{p^2} \operatorname{tr} M_d^2 \\ &= (v-3)r^2 + r + \frac{2v^2}{p^2} \sum\sum_{h<h'}(\delta_{dhh'} - \gamma)^2,\end{aligned}$$

so that $\operatorname{tr} C_d^2$ is minimized by designs that minimize $\sum\sum_{h<h'}(\delta_{dhh'} - \gamma)^2$. Since $\delta_{dhh'}$ is integral, the unique minimum of $\operatorname{tr} C_d^2$ on $0 \leq \gamma < \frac{1}{2}$ is at $\delta_{dhh'} \equiv 0$, achieved only by $\text{ECD}(\bar{\theta})$. For $\gamma = \frac{1}{2}$, any values $\delta_{dhh'} \in \{0, 1\}$ minimize $\operatorname{tr} C_d^2$. $\square$

Now define the F-criterion as the value of the largest eigenvalue of $C_d$ that is not constrained by the setting to equal $r$, that is,

$$\psi_F(C_d) = z_{dr} = r - \frac{ve_{dr}}{p}.$$

Minimizing $\psi_F(C_d)$ over $\mathcal{D}$ is equivalent to maximizing $e_{dr}$. This criterion can be important in establishing Schur-optimality, as shown next.

THEOREM 6. *An $ECD(\bar{\theta})$ is Schur-better than a competitor with a different set of eigenvalues if and only if it is F-equivalent or better than that competitor. Consequently, $ECD(\bar{\theta})$'s are Schur-optimal if and only if they are F-optimal.*

PROOF. Follows from (11) and Lemma 3(i). □

A result of a similar flavor holds for $ECD(\bar{\theta}+1)$'s using the E-criterion. As pointed out by Kunert [15], page 385, designs with eigenvalues $z_{di}$ in the form of Lemma 3(ii) are Schur-best whenever they are E-optimal. For the current problem, this is stated as follows.

THEOREM 7. *An $ECD(\bar{\theta}+1)$ is Schur-better than a competitor with a different set of eigenvalues if and only if it is E-equivalent or better than that competitor. Consequently, $ECD(\bar{\theta}+1)$'s are Schur-optimal if and only if they are E-optimal.*

In establishing necessary and sufficient conditions for Schur-optimality for their respective ECDs in terms of a single eigenvalue, Theorems 6 and 7 provide simple tests for comparing these designs to any other design. An immediate consequence is that an $ECD(\theta)$ with $\theta \notin \{\bar{\theta}, \bar{\theta}+1\}$ is Schur-inferior to at least one of these two competitors. Thus, among ECDs, only these two competitors remain. They are compared to one another in the following corollary.

COROLLARY 8. *$ECD(\bar{\theta})$'s are Schur-better than $ECD(\bar{\theta}+1)$'s if and only if $\gamma \leq \frac{1}{r}$, and $ECD(\bar{\theta}+1)$'s are Schur-better than $ECD(\bar{\theta})$'s if and only if $\gamma \geq \frac{r-1}{r}$.*

PROOF. Simply use (11) and (12) in applying Theorems 6 and 7. □

Though design nonexistence can play a role, Corollary 8 says that one does not expect to find a Schur-optimal design for $\frac{1}{r} \leq \gamma \leq \frac{r-1}{r}$ (Theorem 5 makes for an interesting juxtaposition). Schur-domination can nonetheless be used to eliminate many competitors, as shown next.

3.3. *Schur-inferiority of designs lacking equality of block concurrences.* Given the results of the preceding subsection, the remaining question from a global optimality perspective is if (and when) designs outside the ECD class can be preferable. This question can be effectively pursued by application of Theorems 6 and 7, once workable bounds for the largest and smallest eigenvalues $e_{d1}$ and $e_{dr}$ of the optimality matrix are in place.

LEMMA 9. *Let $d \in \mathcal{D}(v, r; k_1, k_2)$ have concurrence discrepancy matrix $\Delta_d = (\delta_{dhh'})$ and optimality matrix $M_d$. Then:*

(i) $\delta_{d12} \leq 0 \Rightarrow e_{d1} \geq \frac{p}{v} + \gamma - \delta_{d12}$ and $e_{dr} \leq \frac{p}{v} - \gamma + \delta_{d12}$;
(ii) $\delta_{d12} > 0 \Rightarrow e_{d1} \geq \frac{p}{v} - \gamma + \delta_{d12}$ and $e_{dr} \leq \frac{p}{v} + \gamma - \delta_{d12}$.

PROOF. The leading $2 \times 2$ minor of $M_d$, which is $M_{d11} = (\frac{p}{v} + \gamma - \delta_{d12})I - (\gamma - \delta_{d12})J$, has eigenvalues $\frac{p}{v} + \gamma - \delta_{d12}$ and $\frac{p}{v} - \gamma + \delta_{d12}$. A Sturmian separation theorem ([22], page 64) provides the bounds. □

Define design $d$ to be an *almost equal concurrence design*, or AECD, if $\phi_{dhh'} \in \{\bar{\theta}, \bar{\theta} + 1\}$ for all $h \neq h'$ and if each value is attained for some $h, h'$. If any $\phi_{dhh'}$ is not in $\{\bar{\theta}, \bar{\theta} + 1\}$, then $d$ is an *unequal concurrence design*, or UCD. In terms of discrepancies, AECDs have all $\delta_{dhh'} \in \{0, 1\}$, while UCDs have some $\delta_{dhh'} \leq -1$ or $\geq 2$. Depending on $\gamma$, AECDs can be optimal in at least some senses (see Section 4), necessarily ruling out global optimality of ECDs for some $\gamma$. The next few results will show that ECDs often dominate UCDs.

COROLLARY 10. *Suppose $d \in \mathcal{D}(v, b; k_1, k_2)$ is a UCD with $\delta_{dhh'} \leq -\alpha$ for some $1 \leq h \neq h' \leq r$ and some integer $\alpha \geq 1$. Then:*

(i) *$ECD(\bar{\theta})$'s are Schur-better than $d$ if $\gamma \leq \frac{\alpha}{r-2}$;*
(ii) *$ECD(\bar{\theta} + 1)$'s are Schur better than $d$ if $\gamma \geq \frac{r-\alpha-1}{r}$.*

PROOF. For the UCD $d$ as described in the corollary, take $\delta_{d12} \leq -\alpha$. Then from Lemma 9, $e_{d1} \geq \frac{p}{v} + \gamma - \alpha$ and $\frac{p}{v} - \gamma + \alpha \geq e_{dr}$. By Theorem 6, an $ECD(\bar{\theta})$ is Schur-better than $d$ if $\xi_2(\gamma) \geq \frac{p}{v} - \gamma + \alpha \geq e_{dr} \iff \gamma \leq \frac{\alpha}{r-2}$. By Theorem 7, an $ECD(\bar{\theta} + 1)$ is Schur-better than $d$ if $e_{d1} \geq \frac{p}{v} + \gamma - \alpha \geq \xi_2(\gamma - 1) \iff \gamma \geq \frac{r-\alpha-1}{r}$. □

COROLLARY 11. *When $r \leq 4$, all UCDs with some $\delta_{dhh'} \leq -1$ are Schur-inferior to an ECD, and when $r = 5$ or $6$, UCDs with some $\delta_{dhh'} \leq -2$ are Schur-inferior to an ECD.*

The next two corollaries are similarly shown.

COROLLARY 12. *Suppose $d \in \mathcal{D}(v, b; k_1, k_2)$ is a UCD with $\delta_{dhh'} \geq \alpha$ for some $1 \leq h \neq h' \leq r$ and some integer $\alpha \geq 2$. Then:*

(i) *$ECD(\bar{\theta})$'s are Schur-better than $d$ if $\gamma \leq \frac{\alpha}{r}$;*
(ii) *$ECD(\bar{\theta} + 1)$'s are Schur better than $d$ if $\gamma \geq \frac{r-\alpha-1}{r-2}$.*

COROLLARY 13. *When $r \leq 4$, all UCDs with some $\delta_{dhh'} \geq 2$ are Schur-inferior to an ECD, and when $r = 5$ or $6$, UCDs with some $\delta_{dhh'} \geq 3$ are Schur-inferior to an ECD.*

Corollaries 11 and 13 say that optimal (with respect to *any* convex criterion) designs in settings $\mathcal{D}(v, r; k_1, k_2)$ with $r \leq 4$ must be an $\text{ECD}(\bar{\theta})$, an $\text{ECD}(\bar{\theta}+1)$ or an AECD. Optimal designs in settings with $r = 5$ or $6$ must have block concurrence discrepancies $\delta_{dhh'} \in \{-1, 0, 1, 2\}$ for all $1 \leq h \neq h' \leq r$. It is unlikely that such global statements can be much improved. The importance of these results lies in the prevalence of small $r$ in the application of resolvable designs.

**4. E-optimality.** Sufficient conditions, which are necessary given existence, will be developed for E-optimality of designs in $\mathcal{D}(v, r; k_1, k_2)$. The main results, Theorems 14 and 15, will be stated after first defining a subclass where E-optimal designs will be found.

Corollary 16 below will remove the UCDs from E-contention, so that only $\text{ECD}(\bar{\theta})$'s, $\text{ECD}(\bar{\theta}+1)$'s and AECDs need be considered. These designs have all $\delta_{dhh'} \in \{0, 1\}$ and so have $\Delta_d$ which is the adjacency matrix of a simple, undirected graph on $r$ vertices. Any of these designs for which $\Delta_d - J$ is (with suitable ordering of replicates) of the form

$$(13) \quad \begin{pmatrix} -J_{t_1} & & & \\ & -J_{t_2} & & \\ & & \ddots & \\ & & & -J_{t_n} \end{pmatrix}$$

for positive integers $n$ and $t_1 \leq t_2 \leq \cdots \leq t_n$ ($\sum_{i=1}^{n} t_i = r$) is said to be *group-affine*. For group-affine designs, concurrences $\phi_{dhh'}$ are constant ($= \bar{\theta}$) within groups of sizes $t_1, \ldots, t_n$ and are constant ($= \bar{\theta}+1$) between groups. A group-affine design is said to be *uniform* if, for given $n$, the range of group sizes $t_i$ is at most one. For any group-affine design with the number of groups $n \leq n_\gamma = \text{int}(\frac{1}{1-\gamma})$, let $t_{(d)}$ denote the vector of its group sizes $t_i$ arranged in increasing order and with $n_\gamma$ elements, padding with zeros as necessary. For example, if $r = 7$, $\gamma = \frac{7}{9}$ and $d$ has $t_i$ values 1, 3 and 3, then $n_\gamma = 4$ and $t_{(d)} = (0, 1, 3, 3)$. Now the main results can be stated. The phrase "E-Schur-optimal" means "Schur-optimal within the class of all E-optimal designs."

THEOREM 14. *Any group-affine design is E-optimal if $n \leq n_\gamma$ and $\gamma \leq \frac{r-1}{r}$. If any such design exists, then all such designs form the class $\mathcal{D}_E$ of all E-optimal designs. $ECD(\bar{\theta}+1)$'s are Schur-optimal for $\gamma \geq \frac{r-1}{r}$.*

THEOREM 15. *For $d_1, d_2 \in \mathcal{D}_E$, if $t_{(d_1)} \prec t_{(d_2)}$, then $d_1$ is Schur-better than $d_2$. Design $d^* \in \mathcal{D}_E$ is E-Schur-optimal if it is uniform with $n = n_\gamma$.*

Theorem 14 characterizes the class of E-optimal block concurrence structures for $\mathcal{D}(v, r; k_1, k_2)$. Although not much discussed in the literature, the class of E-optimal designs in a given simple block design setting often contains a variety of designs with different information matrices (see, e.g., [18] and [26]). Theorem 14 reveals the same situation for resolvable design settings. This allows other design criteria to be brought to bear: one should choose the best of the E-optimal designs. Theorem 15 does this. In a very strong sense, E-Schur-optimal designs are the best of the E-optimal designs.

The proofs of Theorems 14 and 15 depend on a series of technical results that are developed in the remainder of this section. So as not to overly disrupt the flow of the main line of proof, the longer of the "sub-proofs" are delayed until Appendix A. The first task is to rule out UCDs.

COROLLARY 16. *For all $r \geq 2$ and $0 \leq \gamma < 1$, $ECD(\bar{\theta})$'s are E-better than UCDs.*

PROOF. The maximum eigenvalue of the optimality matrix for an $ECD(\bar{\theta})$ is $\xi_1(\gamma) = \frac{p}{v} + \gamma$. For the UCD $d$, suppose that $\delta_{d12} \leq -\alpha$ for some integer $\alpha \geq 1$. Then by Lemma 9, $e_{d1} \geq \frac{p}{v} + \gamma - \delta_{d12} > \xi_1(\gamma)$ and $ECD(\bar{\theta})$'s are E-better than $d$. If $\delta_{d12} \geq \alpha$ for some integer $\alpha \geq 2$, then another application of the lemma gives $e_{d1} \geq \frac{p}{v} - \gamma + \delta_{d12} > \xi_1(\gamma)$ and, again, $ECD(\bar{\theta})$'s are E-better than $d$. □

Completing the proof of Theorem 14 is a matter of minimizing $e_{d1}$ over all $d$ for which every $\delta_{dhh'} \in \{0, 1\}$. Before providing the details, here is a sketch of what will be done. First, a lower bound for $e_{d1}$ is found in terms of the largest eigenvalue of a specific projection $P\Delta_d P$ of the discrepancy matrix. This bound implies that $ECD(\bar{\theta})$'s are E-superior to all designs for which $P\Delta_d P$ has a positive eigenvalue. The next step is thus to find a necessary and sufficient condition on $\Delta_d$ so that $P\Delta_d P$ is nonpositive definite; this condition turns out to be exactly the definition of group-affine design. Finally, the general form of $M_d$ for group-affine designs is examined in detail to determine which of these competitors are E-optimal, producing the conditions of Theorem 14.

LEMMA 17. *For $d \in \mathcal{D}(v, r; k_1, k_2)$ and with $P = (I - \frac{1}{r}J)$, let $u_{d1}$ and $u_{dr}$ be the maximum and minimum eigenvalues of $P\Delta_d P$, respectively. Then:*

(i) if $u_{d1} > 0$ then $e_{d1} \geq \frac{p}{v} + \gamma + u_{d1}$;
(ii) $e_{dr} \leq \frac{p}{v} + \gamma + u_{dr}$.

The proof of Lemma 17 appears in Appendix A.2. Lemma 17 combined with Theorems 6 and 7 immediately gives the following corollary:

COROLLARY 18. *For $d \in \mathcal{D}(v, r; k_1, k_2)$ and with $P = (I - \frac{1}{r}J)$, let $u_{d1}$ and $u_{dr}$ be the maximum and minimum eigenvalues of $P\Delta_d P$.*

(i) *If $\gamma < -\frac{u_{dr}}{r}$, then $ECD(\bar{\theta})$'s are Schur-better than $d$.*

(ii) *If $u_{d1} > 0$ and $\gamma > (\frac{r - u_{d1} - 1}{r})$, then $ECD(\bar{\theta} + 1)$'s are Schur-better than $d$.*

(iii) *If $u_{d1} > 0$, then $ECD(\bar{\theta})$'s are E-better, but not necessarily Schur-better, than $d$.*

Parts (i) and (ii) of Corollary 18 provide alternative routes (cf. Corollaries 10 and 12) for establishing Schur-domination of ECDs. Part (iii) is the key part of the proof of Theorem 14. If $d$ is not to be eliminated by $ECD(\bar{\theta})$, then $P\Delta_d P$ cannot have a positive eigenvalue, that is, must be nonpositive definite. Recall that the current E-competitors are all $d$ for which $\Delta_d$ is the adjacency matrix of a simple, undirected graph.

LEMMA 19. *Let $A$ be the adjacency matrix for a simple undirected graph of $r$ vertices. $PAP$ is nonpositive definite if and only if $A - J$ may be written (possibly after vertex permutation) in the form (13) for some positive integers $n$ and $t_1 \leq t_2 \leq \cdots \leq t_n$ with $\sum_{i=1}^{n} t_i = r$.*

Combining Lemma 19 with Corollaries 16 and 18(iii), it has now been shown that existence of an $ECD(\bar{\theta})$ implies that an E-optimal design is a group-affine design. Group-affine designs include the $ECD(\bar{\theta})$'s (put $n = 1$, $\Delta_d = 0$) and the $ECD(\bar{\theta} + 1)$ (put $n = r$, $\Delta_d = J - I$) at the extremes for $n$. The proof of Lemma 19 appears in Appendix A.3.

It remains to determine which of the group-affine designs are actually E-best, requiring knowledge of $e_{d1}$ for this class. For any group-affine design, write $M_d = (\frac{p}{v} + \gamma)I + \Delta_d - \gamma J$, where, for given $t_1 \leq \cdots \leq t_n$,

$$
\Delta_d - \gamma J = (-1) \begin{pmatrix} \gamma J_{t_1,t_1} & (\gamma - 1)J_{t_1,t_2} & \cdots & (\gamma - 1)J_{t_1,t_n} \\ (\gamma - 1)J_{t_2,t_1} & \gamma J_{t_2,t_2} & \cdots & (\gamma - 1)J_{t_2,t_n} \\ \vdots & \vdots & \ddots & \vdots \\ (\gamma - 1)J_{t_n,t_1} & (\gamma - 1)J_{t_n,t_2} & \cdots & \gamma J_{t_n,t_n} \end{pmatrix}
$$
(14)
$$\equiv (-1)H_d.$$

An E-optimal design will maximize the minimum eigenvalue of $H_d$. Clearly, $H_d$ has $\sum_{i=1}^{n}(t_i - 1) = r - n$ eigenvalues of zero (corresponding to eigenvectors which are orthogonal contrasts within groups of sizes $t_1, \ldots, t_n$). So all of these designs have $H_d$ with at least one eigenvalue of zero, except $\text{ECD}(\bar{\theta} + 1)$, for which the eigenvalues of $H_d$ are 1 (frequency $r - 1$) and $1 + r(\gamma - 1)$. $\text{ECD}(\bar{\theta} + 1)$ is therefore E-optimal if $1 + r(\gamma - 1) \geq 0$, that is, if $\gamma \geq \frac{r-1}{r}$ (in which case it is Schur-optimal; see Theorem 7). Otherwise, all designs for which $H_d$ is nonnegative definite are E-optimal. Needed now are the eigenvalues of $H_d$ other than zero. Let $D_t$ be the diagonal matrix with diagonal elements $(t_1, \ldots, t_n)$. Lemma 20 is proved in Appendix A.4.

LEMMA 20. *The eigenvalues of $H_d$ specified by* (14) *are* 0 *(with frequency $r - n$) and the eigenvalues of $D_t^{1/2} E D_t^{1/2}$, where $E = I_n - (1 - \gamma)J_n$.*

Consequently, $H_d$ is nonnegative definite if and only if $D_t^{1/2} E D_t^{1/2}$ is nonnegative definite. This is so if and only if $E$ is nonnegative definite, that is, if and only if $1 - \gamma \leq \frac{1}{n}$. Thus, for any $\gamma \leq \frac{r-1}{r}$, the E-optimal designs are all group-affine designs for which the number of groups $n$ is no larger than $\frac{1}{1-\gamma}$, that is, for which $n \leq n_\gamma$. This completes the proof of Theorem 14.

To prove Theorem 15, recall the definition of $t_{(d)}$ given just prior to Theorem 14 and note that now only $d \in \mathcal{D}_E$ is being considered. It is shown in Appendix A.4 that the eigenvalues of $H_d$, aside from $r - n_\gamma$ zeros, are the eigenvalues of $E^{1/2} D_t E^{1/2}$. Thus, the problem is to show that $t_{(d_1)} \prec t_{(d_2)}$ implies the eigenvalues of $E^{1/2} D_{t_{(d_1)}} E^{1/2}$ are majorized by those of $E^{1/2} D_{t_{(d_2)}} E^{1/2}$.

Now, $t_{(d_1)} \prec t_{(d_2)}$ is equivalent to $t_{(d_1)} = S t_{(d_2)}$, where $S$ is a doubly stochastic matrix; $S = \sum_{i=1}^{m} a_i Q_i$ for permutation matrices $Q_i$ and positive numbers $a_i$ summing to 1 ([3], pages 33 and 37). Thus,

$$E^{1/2} D_{t_{(d_1)}} E^{1/2} = E^{1/2} D_{S t_{(d_2)}} E^{1/2}$$

$$= E^{1/2} \left[ \sum_i a_i Q_i D_{t_{(d_2)}} Q_i \right] E^{1/2}$$

$$= \sum_i a_i Q_i [E^{1/2} D_{t_{(d_2)}} E^{1/2}] Q_i,$$

the last equality following because $E^{1/2} = (I - \frac{1}{n}J) + \sqrt{\frac{1}{n} - (1 - \gamma)}J$ commutes with any permutation matrix. This shows that $E^{1/2} D_{t_{(d_1)}} E^{1/2}$ is a symmetrization of $E^{1/2} D_{t_{(d_2)}} E^{1/2}$ and thus the eigenvalues of the first matrix are majorized by those of the second (this follows from [3], page 69). Since $t_{(d^*)} \prec t_{(d)}$ for every $d \in \mathcal{D}_E$, $d^*$ is E-Schur-optimal and Theorem 15 is proved.

**5. Special cases: $(k_1 - k_2) \leq 2$.** In this section, the three important special cases of $k_1$ and $k_2$ being equal or nearly so, $k_2 \in \{k_1, k_1 - 1, k_1 - 2\}$, are investigated in light of the results in Sections 3 and 4. Put $k_2 = k_1 - m$ so that $(k_1 - k_2) \leq 2$ says $m = 0$, 1, or 2, and for any $m$,

$$(15) \qquad \frac{k_1^2}{k_1 + k_2} = \frac{k_1^2}{2k_1 - m} = \frac{k_1}{2} + \frac{m}{4} + \frac{m^2}{4(2k_1 - m)}.$$

Recall that $\bar{\theta}$ is the integer part of (15). The values for $\gamma = \frac{k_1^2}{k_1 + k_2} - \bar{\theta}$ in the corollaries below are easily found using (15).

COROLLARY 21. *For $k_1 = k_2$,*

(i) *if $2|k_1$, then $\gamma = 0$, $(k_1 + k_2)|k_1^2$ and $ECD(\theta^*)$'s are Schur-optimal;*
(ii) *if $2 \nmid k_1$, then $\gamma = \frac{1}{2}$ and $ECD(\bar{\theta})$'s are E-Schur and type-1-optimal.*

When $k_1 = k_2$ and $2|k_1$, Schur-optimality also follows from [1], Corollary 3.4. For $2 \nmid k_1$, the result follows from Theorems 5 and 15.

COROLLARY 22. *For $k_1 - k_2 = 1$,*

(i) *if $2|k_1$, then $\gamma = \frac{v+1}{4v}$ and $ECD(\bar{\theta})$'s are E-Schur and type-1-optimal.*
(ii) *if $2 \nmid k_1$, then $\gamma = \frac{3v+1}{4v}$ and uniform group-affine designs with four groups are E-Schur-optimal if $r \geq 5$. For $r \leq 4$, $ECD(\bar{\theta} + 1)$'s are Schur-optimal.*

As an example, the design in Table 1 is E-Schur-optimal for 9 treatments in 5 replicates with block sizes 4 and 5. The design consisting of the first four replicates is Schur-optimal.

COROLLARY 23. *For $k_1 - k_2 = 2$,*

(i) *if $2|k_1$, then $\gamma = \frac{v+2}{2v}$, and uniform group-affine designs are E-Schur-optimal if the number of groups is 3 for $v = 6$ and 2 for $v \geq 10$. For $r = 2$, $ECD(\bar{\theta} + 1)$'s are Schur-optimal;*
(ii) *if $2 \nmid k_1$, then $\gamma = \frac{1}{v}$, and $ECD(\bar{\theta})$'s are Schur-optimal.*

The Schur-optimality in part (ii) of Corollary 23 follows from Theorem 6 and Lemma 9.

**6. ECDs, balanced arrays and design construction.** Balanced arrays were introduced by Chakravarti [5, 6] as a useful device for fractional factorial designs and they have since been investigated by a plethora of authors, including, of late, Kuriki [16], Fuji-Hara and Miyamoto [9], Sinha, Dhar and Kageyama [25] and Ghosh and Teschmacher [10]. Here, only strength-two arrays on two symbols will be needed. A balanced array of strength 2, $BA(N, m, 2)$, on the symbols 0 and 1, is an $N \times m$ array with the property that for any selection of two columns, the $N$ pairs formed by the rows are $(0,0)$, $(0,1)$, $(1,0)$ and $(1,1)$, with frequencies $\mu_0$, $\mu_1$, $\mu_1$ and $\mu_2$, respectively. Two-symbol orthogonal arrays of strength two are the special case $\mu_0 = \mu_1 = \mu_2$. When the number of 0's in each column is specified (as is the case below), the $\mu_i$'s are all determined by $\theta = \mu_0$ and the array will be denoted $BA(N, m, 2; \theta)$.

Bailey, Monod and Morgan [1] demonstrate the combinatorial equivalence between affine resolvable designs and orthogonal arrays. Their method can be used to express any resolvable design as a combinatorial array, as follows. Given a resolvable design for $v$ treatments in $r$ replicates, each consisting of $s$ blocks, label the blocks within a replicate $0, \ldots, s-1$. Construct a $v \times r$ array by identifying rows of the array with treatments of the resolvable design, and columns with replicates: symbol $j \in \{0, 1, \ldots, s-1\}$ is placed in row $i$, column $h$ if and only if treatment $i$ is in block $j$ of replicate $h$. Theorem 24 is evident.

THEOREM 24. *Each $ECD(\theta)$ in $\mathcal{D}(v, r; k_1, k_2)$ is equivalent to a $BA(v, r, 2; \theta)$.*

A group-affine design with $n$ groups of sizes $t_1, \ldots, t_n$ is a juxtaposition $(BA_1, BA_2, \ldots, BA_n)$ of balanced arrays $BA_i = BA(v, t_i, 2; \theta)$ so that any two columns from different $BA_i$'s form a $BA(v, 2, 2; \theta+1)$. Call such a juxtaposition a *grouped balanced array*, denoted by $GBA(v, (t_1, \ldots, t_n), 2; \theta)$. Constructions for the designs in Corollaries 21–23 of Section 5 are now listed. These are stated in terms of *Hadamard matrices*, that is, orthogonal matrices for which every element is 1 or $-1$. The *Hadamard conjecture* says that a Hadamard matrix exists for every order a multiple of four. Existence is known for all such orders up to 664 and for infinitely many other values (see [8]). A Hadamard matrix is said to be *standardized* if the first row and column are all ones; this can always be achieved.

THEOREM 25. *If $k_1 = k_2$ is even, then a Schur-optimal $ECD(\theta^*)$ corresponds to an $OA(v, r, 2; \frac{v}{2})$. Existence of a Hadamard matrix of order $v$ implies existence of the OA for every $r \leq v - 1$.*

THEOREM 26. *If $k_1 = k_2$ is odd, then a type-1-optimal $ECD(\bar{\theta})$ corresponds to a $BA(v, r, 2; \frac{v-2}{4})$. Existence of a Hadamard matrix of order $v+2$ implies existence of the BA for every $r \leq \frac{v}{2}$.*

PROOF. The value of $\bar{\theta}$ is $\text{int}(\frac{(v/2)^2}{v}) = \frac{v-2}{4}$. Given the standardized Hadamard matrix, permute columns (except the first) so that the second row has 1 in the first $\frac{v+2}{2}$ columns. Delete the first two rows and the first $\frac{v+2}{2}$ columns, then replace $-1$ by 0 throughout. Clearly, the result is $\text{BA}(v, \frac{v}{2}, 2; \frac{v-2}{4})$. □

THEOREM 27. *If $k_1 = k_2 + 1$ is even, then a type-1-optimal $ECD(\bar{\theta})$ corresponds to a $BA(v, r, 2; \frac{v+1}{4})$. Existence of a Hadamard matrix of order $v + 1$ implies existence of the BA for every $r \leq v$.*

PROOF. The value of $\bar{\theta}$ is $\text{int}(\frac{((v+1)/2)^2}{v}) = \frac{v+1}{4}$. Given the standardized Hadamard matrix, delete the first row and column, then replace $-1$ by 0 throughout. The result is $\text{BA}(v, v, 2; \frac{v+1}{4})$. □

THEOREM 28. *If $k_1 = k_2 + 1$ is odd, then a type-1-optimal $ECD(\bar{\theta}+1)$ with $r \leq 4$ corresponds to a $BA(v, r, 2; \frac{v+3}{4})$, which always exists. For $r \geq 5$ and $v \geq 9$, an E-Schur-optimal group-affine design corresponds to a $GBA(v, (t_1, t_2, t_3, t_4), 2; \frac{v-1}{4})$ with $t_4 - t_1 \leq 1$. Existence of a Hadamard matrix of order $v + 3$ with a $4 \times r$ submatrix of the form*

$$\begin{pmatrix} 1_{t_1} & 1_{t_2} & 1_{t_3} & -1_{t_4} \\ 1_{t_1} & 1_{t_2} & -1_{t_3} & 1_{t_4} \\ 1_{t_1} & -1_{t_2} & 1_{t_3} & 1_{t_4} \\ -1_{t_1} & 1_{t_2} & 1_{t_3} & 1_{t_4} \end{pmatrix}$$

*implies existence of the GBA.*

PROOF. The value of $\bar{\theta}$ is $\text{int}(\frac{((v+1)/2)^2}{v}) = \frac{v-1}{4}$. For $r \leq 4$, take $OA(v-1, r, 2)$ on $\{0, 1\}$ and add one row of 0's to get $\text{BA}(v, r, 2; \frac{v+3}{4})$. Given the assumed Hadamard matrix, delete $v + 3 - r$ columns not containing the submatrix, then delete the submatrix and add a row of 1's. Replacing $-1$ by 0 throughout gives the GBA. □

The maximum $r$ admitted for a given $v$ by the Hadamard construction in Theorem 28 depends on the particular Hadamard matrix chosen. All nonisomorphic Hadamard matrices are known up through order 28 and, for these, a complete search has produced $(v, r) = (9, 5), (13, 5), (17, 9), (21, 9), (25, 13)$; the first four of these appear in Appendix B. A search of a few known Hadamard matrices (compiled by N. J. A. Sloane at www.research.att.com/˜njas/hadamard/) of orders up to 48 has further produced $(v, r) = (29, 13), (33, 16), (37, 18), (41, 24), (45, 19)$. Clearly there is room for more work to be done here.

THEOREM 29. *If $k_1 = k_2 + 2$ is even, for $v \geq 10$ an E-Schur-optimal group-affine design corresponds to a $GBA(v, (t_1, t_2), 2; \frac{v+2}{4})$ with $t_2 - t_1 \leq 1$. Existence of a Hadamard matrix of order $v + 2$ implies existence of the GBA for every $t_1 + t_2 = r \leq \frac{v}{2} + 1$.*

PROOF. The value of $\bar{\theta}$ is $\text{int}(\frac{((v+2)/2)^2}{v}) = \frac{v+2}{4}$. Given the standardized Hadamard matrix of order $v + 2 = 4h$, permute columns (except the first) so that the first three rows are

$$\begin{pmatrix} 1_h & 1_h & 1_h & 1_h \\ 1_h & -1_h & 1_h & -1_h \\ 1_h & -1_h & -1_h & 1_h \end{pmatrix}.$$

Now delete the first $2h$ columns, delete the first three rows, then add one row of $-1$'s and finally replace $-1$ by 0 throughout. The result is $\text{GBA}(v, (\frac{v+2}{4}, \frac{v+2}{4}), 2; \frac{v+2}{4})$ with two juxtaposed $\text{BA}(v, \frac{v+2}{4}, 2; \frac{v+2}{4})$ being the first $h$ and last $h$ columns. For smaller $r$, delete $\text{int}(h - \frac{r}{2})$ columns from the first component BA and $\text{int}(h - \frac{r+1}{2})$ from the second. $\square$

THEOREM 30. *If $k_1 = k_2 + 2$ is odd, a Schur-optimal $ECD(\bar{\theta})$ corresponds to a $BA(v, r, 2; \frac{v+4}{4})$. Existence of a Hadamard matrix of order $v$ implies existence of the BA for every $r \leq v - 1$.*

PROOF. The value of $\bar{\theta}$ is $\text{int}(\frac{((v+2)/2)^2}{v}) = \frac{v+4}{4}$. Given the standardized Hadamard matrix, delete the first row and column, add a row of $-1$'s and then replace $-1$ by 0 throughout. The result is $\text{BA}(v, v-1, 2; \frac{v+4}{4})$. $\square$

**7. Comments.** By exploiting properties of the dual, optimality theory for resolvable designs with two blocks per replicate has been developed. Section 3 establishes conditions for Schur-optimality of equal concurrence designs depending on the block discordancy coefficient $\gamma$. For $\gamma$ where Schur-optimality is not established, the class of competitors has been significantly narrowed via Schur-ordering for small $r$. Section 4 characterizes the class of all E-optimal designs whenever any group-affine design exists and further determines the Schur-best of the E-optimal designs. Sections 5 and 6 apply these results for the important cases $k_1 - k_2 \leq 2$, including explicit design constructions employing an equivalence with balanced arrays. It is evident from the constructions that many other designs, eliminated by the Schur-domination argument in Theorem 15, do exist.

For large $(r \geq v)$ replication, the problem can be quite different and optimality work would proceed based on treatment, rather than block, concurrences. Here is a simple construction for that case: given a BIBD for $v$ treatments in $b$ blocks of size $k$, there is a Schur-optimal resolvable design in $\mathcal{D}(v, b; k, v - k)$ formed by the blocks of the BIBD and their complements. While it is possible that these designs are equivalent to balanced arrays (when the starting BIBD is symmetric), it is also true that the resulting $\phi_{dhh'}$ can be widely dispersed, for they are determined by the block concurrence counts for the underlying BIBD.

The optimality route is clear for $r \leq v - 1$ and $\gamma \leq \frac{1}{2}$, but for larger $\gamma$ there are still unanswered questions. For instance, A- and E-optimality need not coincide and the problem of determining an A-best design remains open. For $r \leq 4$, the authors have solved the A-optimality problem in its entirety, including construction, and for larger $r$, have done this for the special cases of $k_1 - k_2 \leq 2$. These and related results will be reported elsewhere.

## APPENDIX A: PROOFS

**A.1. Proof of Lemma 1.** Let $D_s^{1/2}$ be the diagonal matrix of square roots of block sizes and write $B_d = N_d D_s^{-1/2}$. Multiplying $C_d$ by $\frac{1}{r}$ and right and left multiplying $C_{\text{dual}}$ by $D_s^{-1/2}$, equations (1) and (2) become

$$(16) \qquad \frac{1}{r} C_d = I - \frac{1}{r} N_d D_s^{-1} N_d' = I - \frac{1}{r} B_d B_d' = C_d^*$$

and

$$(17) \quad D_s^{-1/2} C_{\text{dual}} D_s^{-1/2} = I - \frac{1}{r} D_s^{-1/2} N_d' N_d D_s^{-1/2} = I - \frac{1}{r} B_d' B_d = C_{\text{dual}}^*.$$

Since (1) and (16) differ only by a constant and since the nonzero eigenvalues of $B_d B_d'$ and $B_d' B_d$ are identical, the eigenvalues of $C_d$ can be found from those of $C_{\text{dual}}^*$ in (17), which depends on block concurrences and block sizes.

Continuing, write $C_{\text{dual}}^* = I - \frac{1}{r} A_d$ for $A_d = B_d' B_d = D_s^{-1/2} N_d' N_d D_s^{-1/2}$. Regardless of the design $d$, $A_d D_s^{1/2} 1 = r D_s^{1/2} 1$; $r$ is an eigenvalue of $A_d$ corresponding to the zero eigenvalue common to $C_{\text{dual}}$ and $C_d$. One term of the spectral decomposition of $A_d$ is then

$$(18) \qquad \frac{r(D_s^{1/2} 1)(D_s^{1/2} 1)'}{(D_s^{1/2} 1)'(D_s^{1/2} 1)} = \frac{1}{(k_1 + k_2)} \left[ J \otimes \begin{pmatrix} k_1 & \sqrt{k_1 k_2} \\ \sqrt{k_1 k_2} & k_2 \end{pmatrix} \right],$$

where $J$ is an $r \times r$ matrix of ones. Subtract (18) from $A_d$ and denote the result $A_d^*$. Then a bit of manipulation, employing the fact that all four block concurrence counts for blocks in replicate $h$ with blocks in replicate $h'$ are determined by $\phi_{dhh'}$, gives

$$(19) \qquad A_d^* = \frac{1}{k_1 k_2} M_d \otimes \begin{pmatrix} k_2 & -\sqrt{k_1 k_2} \\ -\sqrt{k_1 k_2} & k_1 \end{pmatrix}.$$

Since the eigenvalues of the $2 \times 2$ matrix in (19) are 0 and $k_1 + k_2$, the $b = 2r$ eigenvalues of $A_d^*$ are $r$ copies of 0 and $\frac{v}{k_1 k_2}$ times the $r$ eigenvalues of $M_d$. The eigenvalues of $C_d$ as stated in the lemma are now immediate.

**A.2. Proof of Lemma 17.** The lower bound on $e_{d1}$ follows from

$$e_{d1} = \max_{x'x=1} x'M_d x$$

$$= \max_{x'x=1} x'\left[\left(\frac{p}{v}+\gamma\right)I - \gamma J + \Delta_d\right]x$$

$$\geq \max_{\substack{x'x=1 \\ x'1=0}} x'\left[\left(\frac{p}{v}+\gamma\right)I - \gamma J + \Delta_d\right]x$$

$$= \frac{p}{v} + \gamma + \max_{\substack{x'x=1 \\ x'1=0}} x'\Delta_d x$$

$$= \frac{p}{v} + \gamma + \max_{\substack{x'x=1 \\ x'1=0}} x'P\Delta_d P x$$

$$= \frac{p}{v} + \gamma + \max_{x'x=1} x'P\Delta_d P x$$

$$= \frac{p}{v} + \gamma + u_{d1}.$$

The penultimate equality holds since $u_{d1} > 0$ and 1 is an eigenvector of $P\Delta_d P$ with eigenvalue 0. Likewise,

$$e_{dr} = \min_{x'x=1} x'M_d x$$

$$= \min_{x'x=1} x'\left[\left(\frac{p}{v}+\gamma\right)I - \gamma J + \Delta_d\right]x$$

$$\leq \frac{p}{v} + \gamma + \min_{\substack{x'x=1 \\ x'1=0}} x'P'\Delta_d P x$$

$$= \frac{p}{v} + \gamma + u_{dr}.$$

The last equality holds provided $u_{dr} < 0$, for reasons similar to those above. If $u_{dr} > 0$, the bound still holds, since $e_{dr} \leq \frac{\text{tr}(M_d)}{r} = \frac{p}{v} \leq \frac{p}{v} + \gamma + u_{dr}$.

**A.3. Proof of Lemma 19.** Suppose that $A - J$ is of the suggested form. Write $t_{(i)}$ for $t_1 + t_2 + \cdots + t_i$ and $t_{(0)} = 0$. Then $PAP = P(A-J)P$, implying that

$$\max \text{eig}(PAP) = \max_{x'x=1} x'P(A-J)Px$$

$$= \max_{y=Px, x'x=1} y'(A-J)y$$

$$= \max_{y=Px, x'x=1} -\sum_{i=1}^{t}\left(\sum_{j=t_{(i-1)}+1}^{t_{(i)}} y_j\right)^2,$$

which is clearly nonpositive.

Now, subscripting by the dimension, suppose that $P_r A_r P_r$ is nonpositive definite (npd). Exhaustive enumeration shows that $A_r - J_r$ must have the form (13) for $r = 3, 4, 5$. Assuming nonpositivity implies that the form must hold for a given $r$, it will be shown that the same implication holds for $r+1$. Denoting the upper left-hand side $r \times r$ submatrix of $A_{r+1}$ by $A_r$ (which is itself an adjacency matrix), it is claimed that $P_{r_1} A_{r+1} P_{r+1}$ is npd $\Rightarrow$ $P_r A_r P_r$ is npd. If not, then there exists $x$ such that $x' P_r A_r P_r x > 0$. Write $y' = (x' P_r, 0)$. Then $y'1 = 0 \Rightarrow y' P_{r+1} A_{r+1} P_{r+1} y = y' A_{r+1} y = x' A_r x > 0$, a contradiction.

Since $P_r A_r P_r$ is npd, the induction hypothesis says that

$$A_{r+1} - J_{r+1} = \begin{pmatrix} -J_{t_1} & & & & \\ & -J_{t_2} & & & a \\ & & \ddots & & \\ & & & -J_{t_n} & \\ \hline & & a' & & -1 \end{pmatrix}$$

for some vector $a_{r \times 1}$ of 0's and $-1$'s. Indeed, by the induction hypothesis, every $s \times s$ principal minor of $A_{r+1} - J_{r+1}$ must have the form (13), so $a$ may be partitioned as $a' = (a'_1, a'_2, \ldots, a'_n)$, where $a_i$ is either $-1_{t_i}$ or $0_{t_i}$. If $n = 1$, the proof is done. If $n = 2$ and $a_2 = -1_{t_2}$, then

$$A_{r+1} - J_{r+1} = \begin{pmatrix} -J_{t_1} & 0 & a_1 \\ 0 & -J_{t_2} & -1_{t_2} \\ a'_1 & -1'_{t_2} & -1 \end{pmatrix}.$$

If now $a_1 = -1_{t_1}$, then $A_{r+1} - J_{r+1}$ contains a principal minor,

$$M_3 = (-1) \begin{pmatrix} 1 & 0 & 1 \\ 0 & 1 & 1 \\ 1 & 1 & 1 \end{pmatrix},$$

which is not of the form (13), contradicting the fact that the result holds for $r = 3$. Thus, at most one of $a_1, a_2$ is nonzero. It follows immediately that for $n > 2$, the same statement holds for any pair $a_i, a_{i'}$. Permuting rows and columns as needed, it can be assumed that $a_i = 0_{t_i}$ for $i < n$ and consequently that $A_{r+1} - J_{r+1}$ has the form (13).

**A.4. Eigenvalue equations for $H_d$ and proof of Lemma 20.** Any vector of the form $(c'_1, c'_2, \ldots, c'_n)$, where $c_i \in \Re^{t_i}$ is either a contrast vector or the zero vector, is an eigenvector of $H_d$ with eigenvalue 0. Consequently, all other eigenvectors are of the form $e' = (x_1 1'_{t_1}, x_2 1'_{t_2}, \ldots, x_n 1'_{t_n})$ for some scalars $x_1, \ldots, x_n$. The first equation in the system $H_d e = \lambda e$ is $\gamma t_1 x_1 - (1-\gamma) t_2 x_2 - \cdots - (1-\gamma) t_n x_n = \lambda x_1$; the other equations may be written similarly in order to see that $H_d e = \lambda e$ are equivalent to $[D_t - (1-\gamma) 1_n t'] x = \lambda x$, where $t = (t_1, \ldots, t_n)'$ and $x = (x_1, \ldots, x_n)'$. Thus, the remaining eigenvalues of $H_d$ are the right eigenvalues of $D_t - (1-\gamma)1_n t'$. Now, $D_t$ is positive definite and $E$ is nonnegative definite for $\gamma \geq \frac{n-1}{n}$, so both have symmetric square root matrices and $D_t$ is invertible. Thus,

$$\begin{aligned}
|D_t - (1-\gamma)1_n t' - \lambda I| &= 0 \\
\iff \quad |D_t - (1-\gamma)1_n t' - \lambda I||D_t^{-1}| &= 0 \\
\iff \quad |E - \lambda D_t^{-1}| &= 0 \\
\iff \quad |D_t^{1/2} E D_t^{1/2} - \lambda I| &= 0 \quad \text{(proving Lemma 20)} \\
\iff \quad |E^{1/2} D_t E^{1/2} - \lambda I| &= 0 \quad \text{(needed in Theorem 15)},
\end{aligned}$$

the last step following because $GG'$ and $G'G$ have the same eigenvalues for any square $G$. Now, in the proof of Theorem 15 some elements of $t$ are allowed to be zero (without loss of generality, $t_1 = \cdots = t_z = 0$ for integer $z \geq 1$), in which case $D_t$ is not invertible. In this case there are $z$ additional zero eigenvalues plus a reduced system of $n - z$ equations in $t_{z+1}, \ldots, t_n$. It is easy to see that $|D_t^{1/2} E D_t^{1/2} - \lambda I| = |0_{z,z} - \lambda I_z||\tilde{D}_t^{1/2} \tilde{E} \tilde{D}_t^{1/2}|$, where $\tilde{E}$ and $\tilde{D}_t$ are the lower-right submatrices of $E$ and $D_t$ of order $n - z$. Thus, the $n$ eigenvalues sought are still those of $D_t^{1/2} E D_t^{1/2}$ and thus of $E^{1/2} D_t E^{1/2}$.

## APPENDIX B: UNIFORM GBAS WITH FOUR GROUPS

GBA(9, (1, 1, 1, 2), 2; 2)  GBA(13, (1, 1, 1, 2), 2; 3)

```
1 1 1 1 1         1 1 1 1 1
0 0 0 0 0         1 1 1 1 1
0 1 0 0 1         0 1 1 0 0
0 1 0 1 0         0 0 0 1 0
0 0 1 1 0         1 0 0 1 0
1 0 0 1 0         0 0 0 0 1
0 0 1 0 1         1 0 0 0 1
1 1 1 0 0         0 0 1 1 0
1 0 0 0 1         1 0 1 0 0
                  0 0 1 0 1
                  0 1 0 1 0
                  0 1 0 0 1
                  1 1 0 0 0
```

$$\text{GBA}(17, (2, 2, 2, 3), 2; 4) \qquad \text{GBA}(21, (2, 2, 2, 3), 2; 5)$$

```
1 1 1 1 1 1 1 1 1      1 1 1 1 1 1 1 1 1
0 0 0 0 0 0 0 0 0      1 1 1 1 1 1 1 1 1
1 0 0 1 0 0 1 0 1      1 1 0 0 0 0 1 1 0
0 0 1 1 0 1 1 0 0      1 0 1 0 1 0 1 0 0
1 0 0 0 1 1 0 0 1      0 1 1 0 1 0 0 0 1
1 0 1 1 1 0 0 1 0      0 0 0 1 0 1 1 0 0
0 1 0 0 1 0 1 1 0      0 0 1 0 1 0 0 1 0
0 1 0 1 1 1 1 0 0      1 0 0 0 0 0 0 1 1
1 0 0 0 0 1 1 1 0      0 1 0 1 0 1 0 1 0
1 1 1 0 0 1 0 0 1      0 1 0 0 0 0 1 0 1
1 0 1 0 1 0 1 0 0      1 0 0 1 0 1 0 0 1
0 0 0 1 0 1 0 1 1      0 1 1 0 0 1 0 0 1
0 1 1 0 0 1 0 1 0      1 0 1 0 0 1 0 0 1
1 1 0 1 0 0 0 1 0      0 1 0 0 1 1 0 1 0
0 1 0 1 1 0 0 0 1      0 1 1 1 0 0 1 0 0
0 0 1 0 1 0 0 1 1      0 0 1 0 0 1 1 1 0
0 1 1 0 0 0 1 0 1      1 0 1 1 0 0 0 1 0
                       1 0 0 0 1 1 1 0 0
                       0 0 0 1 1 0 1 0 1
                       1 1 0 1 1 0 0 0 0
                       0 0 0 1 1 0 0 1 1
```

**Acknowledgments.** This research is partially based on work in Reck's doctoral dissertation [23], for which Morgan was director. Morgan completed his Ph.D. under the direction of I. M. Chakravarti. In dedicating this paper to the memory of Professor Chakravarti, we humbly acknowledge the ongoing impact of his work on our own and on that of the statistical and mathematical communities at large. We wish to thank the referees for comments leading to decided improvement of the presentation.

Department of Statistics  
Virginia Tech  
Blacksburg, Virginia 24061-0439  
USA  
E-mail: jpmorgan@vt.edu

Genetic Analysis  
GlaxoSmithKline  
Research Triangle Park, North Carolina 27709  
USA  
E-mail: brian.h.reck@gsk.com